\documentclass{article}
\usepackage{amsmath, amsthm}
\usepackage{amssymb}
\usepackage{textcomp}
\newtheoremstyle{theorem}
  {10pt}          
  {10pt}  
  {\sl}  
  {\parindent}     
  {\bf}  
  {. }    
  { }    
  {}     
\theoremstyle{theorem}
\newtheorem{theorem}{Theorem}[section]

\newtheorem{lemma}[theorem]{Lemma}

\newtheorem{example}[theorem]{Example}
\newtheoremstyle{defi}
  {10pt}          
  {10pt}  
  {\rm}  
  {\parindent}     
  {\bf}  
  {. }    
  { }    
  {}     
\theoremstyle{defi}
\newtheorem{definition}[theorem]{Definition}

\begin{document}

\begin{center}
\large{ Variational Optimal-Control Problems with Delayed
Arguments on Time Scales }
\end{center}
\begin{center}
 \textbf{Thabet Abdeljawad (Maraaba), }\textbf{ Fahd Jarad, \\}
\textbf{Dumitru Baleanu}\footnote{On leave of absence from
Institute of Space Sciences, P.O.BOX, MG-23, R
76900,Magurele-Bucharest, Romania, Emails: dumitru@cankaya.edu.tr,
baleanu@venus.nipne.ro}\\
 Department of Mathematics and Computer Science\\
  \c{C}ankaya University, 06530 Ankara, Turkey\\ \vspace{0.5cm}
\end{center}

\begin{abstract}
This article deals with variational optimal-control problems on
time scales in the presence of delay in the state variables. The
problem is considered on a time scale unifying the discrete, the
continuous and the quantum cases. Two examples  in the discrete
and quantum cases are analyzed to illustrate our results.

\end{abstract}


\section{Introduction} \label{s:1}

The calculus of variations interacts deeply with some branches of
sciences and engineering e.g. geometry, economics, electrical
engineering and so on \cite{Rose}. Optimal control problems appear
in  various disciplines of sciences and engineering as well
\cite{Young}.

Time scale calculus was initiated by Hilger  ( see
ref. \cite{Hilger} and the references therein) having in mind   to
unify two existing approaches of dynamic models-difference and
differential equations into a general framework. This kind of
calculus  can be used to model dynamic processes whose time
domains are more complex than the set of integers  or real numbers
\cite{Martinbook}. Several potential applications for this  new
theory were reported (see for example
Refs.\cite{Martinbook}, \cite{Atici}, \cite{Guseinov} and the
references therein). Many researchers studied calculus of variations on time scales. Some of them followed the delta approach and some others followed the nabla approach (see for example Refs. \cite{Torreshigher}, \cite{Torresisoperimetric}, \cite{Mar}, \cite{Torresremarks}, \cite{Torresstrong} and \cite{TorresNoeth} ).

 It is well known that the presence of delay is of great importance in
applications. For example, its appearance in
dynamic equations, variational problems and optimal control
problems may affect the stability of solutions. Very recently,
some authors payed the attention to the importance of imposing the
delay in fractional variational problems  \cite{dth}. The
non-locality of the fractional operators and the presence of delay
as well may give better results for problems involving the
dynamics of complex systems. To the best of our knowledge, there
is no work in the direction of variational optimal-control
problems with delayed arguments on time scales.

Our aim in this article is to obtain the Euler-Lagrange equations
for a functional, where  the state variables of its Lagrangian are
defined on a time scale whose backward jumping operator is
$\rho(t)=qt-h, ~q>0,~ h\geq0$. This time scale, of course, absorbs
the discrete, the continuous and the quantum cases. The state variables
of this Lagrangian  allow the presence of delay as well. Then, we
generalize the results to the n-dimensional case. Dealing with
such a very general problem  enables us to recover many previously
obtained results \cite{Bliss, Cad, Gas, Nat}.

The structure of the article is as follows:

In section 2 basic definitions and preliminary concepts about time
scale are presented. The nabla time scale derivative analysis is
followed there. In section 3 the Euler-Lagrange equations into one
unknown function and then in the n-dimensional case are obtained.
In section 4 the variational optimal control problem is proposed
and solved. In section 5 the results obtained in the previous
sections are particulary studied in the discrete and quantum
cases, where two examples are analyzed in details.  Finally,
section 6 contains our conclusions.

\section{Preliminaries}\label{s2}
A time scale is an arbitrary closed subset of the real line
$\mathbb{R}$. Thus the real numbers and the natural numbers,
$\mathbb{N}$, are examples of a time scale. Throughout this article,
and following \cite{Martinbook}, the time scale will be denoted by
$\mathbb{T}$.
 The forward jump
operator $\sigma:\mathbb{T} \rightarrow \mathbb{T}$ is defined by
$$\sigma(t):= \inf\{s\in \mathbb{T}:s>t\},$$
while the backward jump operator $\rho:\mathbb{T}\rightarrow
\mathbb{T}$ is defined by
$$\rho(t):=\sup\{s\in \mathbb{T}: s<t\},$$
where, $\inf \emptyset= \sup \mathbb{T}$ (i.e. $\sigma(t)= t$ if
$\mathbb{T}$ has a maximum $t$) and $\sup \emptyset = \inf \mathbb{T}$
(i.e $\rho(t)= t$ if $\mathbb{T}$ has a minimum $t$). A point
$t\in \mathbb{T}$ is called right-scattered if $t< \sigma(t)$,
left-scattered if $\rho(t)< t$ and isolated if $\rho(t)< t <
\sigma (t)$. In connection we
define the backward graininess function $\nu:
\mathbb{T}\rightarrow [0,\infty)$ by
$$\nu(t)=t-\rho(t).$$
In order to define the backward time scale derivative down, we
need the set $\mathbb{T}_\kappa$ which is derived from the time
scale $\mathbb{T}$ as follows: if $\mathbb{T}$ has a
right-scattered minimum $m$, then $\mathbb{T}_\kappa=
\mathbb{T}-\{m\}$. Otherwise, $\mathbb{T}_\kappa =\mathbb{T}$.

\begin{definition}\cite{AGus} Assume $f:\mathbb{T}\rightarrow \mathbb{R}$ is a
function and $t
\in~ \mathbb{T}_\kappa $. Then the backward time-scale
derivative $f^\nabla (t)$ is the number (provided it exists) with
the property that given any $\epsilon > 0,$ there exists a
neighborhood $U$ of $t$ (i.e, $U=(t-\delta,t+\delta)$ for some
$\delta > 0$) such that
\begin{equation} \label{s1}
|[f(s)-f(\rho(t))]-[s-\rho(t)]|\leq \epsilon |s-\rho(t)|
~\mbox{for~all}~s \in U
\end{equation}
Moreover, we say that $f$ is (nabla) differentiable on $\mathbb{T}_\kappa$ provided that $f^\nabla (t)$ exists for all $t \in
\mathbb{T}_\kappa$.
\end{definition}
The following theorem is Theorem 3.2 in \cite{Advances} and  an analogue to Theorem 1.16 in
\cite{Martinbook}.

\begin{theorem} \label{boh} \cite{AGus}
Assume $f:\mathbb{T}\rightarrow \mathbb{R}$ is a function and $t
\in \mathbb{T}_\kappa $. Then we have the following:

(i) If $f$ is differentiable at $t$ then $f$ is continuous at $t$.

(ii)If $f$ is continuous at $t$ and $t$ is left-scattered, then
$f$ is differentiable at $t$ with
$$f^\nabla(t)=\frac{f(t)-f(\rho(t))}{\nu(t)}$$

(iii) If $t$ is left-dense, then f is differentiable at $t$ if and
only if the limit
$$\lim_{s\rightarrow t}\frac{f(t)-f(s)}{t-s}$$ exists as a finite number.
In this case
$$f^\nabla(t)=\lim_{s\rightarrow t}\frac{f(t)-f(s)}{t-s}$$
(iv) If $f$ is $\nabla$- differentiable at $t$, then
$$f(t)=f(\rho(t))+\nu(t)f^\nabla(t)$$
\end{theorem}
\begin{example} \label{one}
(i) $\mathbb{T}=\mathbb{R}$ or any any closed interval (The
continuous case) $\sigma(t)=\rho(t)=t,~\nu(t)=0$ and
$f^\nabla(t)=f^\prime(t)$.

(ii) $\mathbb{T}=h\mathbb{Z},~h > 0$ or any subset of it. (The
difference calculus, a discrete case) $\sigma(t)=t+h,~
\rho(t)=t-h,~\nu(t)=h$ and $f^\nabla(t)=\nabla_hf(t)=f(t)-f(t-h)$.

(iii) $\mathbb{T}=\mathbb{T}_q=\{q^n:n \in \mathbb{Z}\} \cup
\{0\}$, $0<q<1$, (quantum calculus)

$\sigma(t)=q^{-1}t, \rho(t)=qt,~\nu(t)=(1-q)t$ and
$f^\nabla(t)=\nabla_q f(t)=\frac{f(t)-f(qt)}{(1-q)t}$.

(iv)  $\mathbb{T}=\mathbb{T}_q^h=\{q^k-\sum_{i=0}^{k-2}q^ih:k \geq
2,k \in \mathbb{N}\}\cup \{\frac{-h}{1-q}\},~~0< q<1,~h>0$
(unifying the difference calculus  and quantum calculus).
$\sigma(t)=q^{-1}(t+h), \rho(t)=qt-h,~\nu(t)=(1-q)t+h$ and
$f^\nabla(t)=\nabla_q^h f(t)=\frac{f(t)-f(qt-h)}{(1-q)t+h}$. If
$\alpha_0 \in \mathbb{N}$ then
$\rho^{\alpha_0}(t)=q^{\alpha_0}t-\sum_{k=0}^{\alpha_0-1}q^k h$
and so $\nabla_q^h \rho^{\alpha_0}(t)=q^{\alpha_0}$. Note that in
this example the backward operator is of the form $\rho(t)=ct+d$
and hence $\mathbb{T}_q^h$ is an element of the class $H$ of time
scales that contains the discrete, the usual and the quantum
calculus (see \cite{Nat}).
\end{example}

\begin{theorem} \label{t} Suppose $f,g: \mathbb{T}\rightarrow \mathbb{R}$
are nabla
differentiable at $t\in ~\mathbb{T}_\kappa$. Then,

1. the sum $f+g:\mathbb{T}\rightarrow \mathbb{R}$ is nabla
differentiable at $t$ and
$(f+g)^\nabla(t)=f^\nabla(t)+g^\nabla(t);$

2. for any $\lambda \in \mathbb{R}$, the function $\lambda
f:\mathbb{T}\rightarrow \mathbb{R}$ is nabla differentiable at $t$
and $(\lambda f)^\nabla(t)= \lambda f^\nabla (t)$;

3. the product $fg:\mathbb{T}\rightarrow \mathbb{R}$ is nabla
differentiable at $t$ and
$$(fg)^\nabla=f^\nabla(t)g(t)+f^\rho(t)g^\nabla(t)=f^\nabla(t)g^\rho(t)+f(t)g^\nabla(t)$$
\end{theorem}
For the proof of the following lemma we refer to \cite{thabet}:
\begin{lemma} \label{th}
Let $\mathbb{T}$ be an $H-$time scale (In particular
$\mathbb{T}=\mathbb{T}_q^h$), $f:\mathbb{T}\rightarrow \mathbb{R}$
two times nabla differentiable function and
$g(t)=\rho^{\alpha_0}(t),$ for $\alpha_0 \in \mathbb{N}$. Then
$$(f \circ g)^\nabla (t)= f^\nabla(g(t)).g^\nabla (t),~~t
\in\mathbb{T}_\kappa$$
\end{lemma}
Throughout this article we use for the time scale derivatives and
integrals the symbol $\nabla_q^h$ which is inherited from the time
scale $\mathbb{T}_q^h$. However, our results are true also for the
$H-$ time scales (those time scales whose jumping operators have
the form $at+b$). The time scale $\mathbb{T}_q^h$ is a natural
example of an $H-$ time scale.
\begin{definition} \label{d1}
A function $F:\mathbb{T}\rightarrow \mathbb{R}$ is called a nabla
antiderivative of $f:T\rightarrow \mathbb{R}$ provided
$F(t)=f(t),$ for all $t\in \mathbb{T}_\kappa$. In this case,
for $a,b\in \mathbb{T}$,  we write $$\int_a^bf(t)\nabla t :=
F(b)-F(a)$$
\end{definition}

The following lemma which extends  the fundamental lemma of
variational analysis on time scales with nabla derivative is
crucial in proving the main results.
\begin{lemma} \label{ftcv}
Let $g \in C_{ld}, g:[a,b]\rightarrow \mathbb{R}^n$. Then
\begin{equation} \label{f1}
\int_a^b g^T(t)\eta^\nabla (t)\nabla t ~~\mbox{for all}~ \eta \in
C_{ld}^1~~\mbox{with}~~ \eta (a)=\eta(b)=0
\end{equation}
holds if and only if
\begin{equation} \label{f2}
g(t)\equiv c ~~ \mbox{on}~~[a,b]_\kappa~~\mbox{for some}~~c \in
\mathbb{R}^n
\end{equation}
\end{lemma}
The proof can be achieved by following as in the proof of Lemma
4.1 in \cite{Mar}, ( see also \cite{Nat}).
\section{First order Euler-Lagrange equation with delay}\label{s3}
We consider the $\mathbb{T}_q^h$-integral functional
$J:S\rightarrow \mathbb{R}$,
\begin{equation}\label{m1}
J(y)=\int_a^b L(x,y^\rho(x),\nabla_q
^hy(x),y^\rho(\rho^{\alpha_0}(x)),\nabla_q^hy((\rho^{\alpha_0}(x))\nabla_q^hx
\end{equation}
where
$$a,b\in \mathbb{T}_q^h,a < \rho^{\alpha_0}(b)< b $$

$$L:[a,b]\times (\mathbb{R}^n)^4\rightarrow \mathbb{R}~\mbox{and}
~y^\rho(x)=y(\rho(x))$$

and

$$S=\{y:[\rho^{\alpha_0}(a),b]\rightarrow
\mathbb{R}^n:y(x)=\varphi(x)~(\forall x \in [\rho^{\alpha_0}(a),
a])~\texttt{and}~y(b)=c_0\}$$

We shall shortly write : $$L(x)\equiv L(x,y^\rho(x),\nabla_q^h
y(x),y^\rho(\rho^{\alpha_0}(x)),\nabla_q^hy((\rho^{\alpha_0}(x))$$

We calculate the first variation of the functional $J$ on the
linear manifold $S$: Let $\eta \in
H=\{h:[\rho^{\alpha_0}(a),b]\rightarrow \mathbb{R}^n: h(x)=0
~(\forall x \in[\rho^{\alpha_0}(a),a]\cup \{b\} )~\}$, then
$$\delta J(y(x),\eta(x))=\frac{d}{d\epsilon}J(y(x)+\epsilon
\eta(x))|_{\epsilon =0}$$
\begin{equation} \label{cal1}
\int_a^b [\partial_1 L(x) \eta^\rho(x)+\partial_2 L (x) \nabla_q^h \eta (x)+ \partial_3 L(x)
\eta^\rho(\rho^{\alpha_0}(x))+q^{\alpha_0} \partial_4 L(x) \nabla_q^h
\eta(\rho^{\alpha_0}(x)) ]\nabla_q^hx
\end{equation}
where
$$\partial_1L=\frac{\partial L}{\partial(y^\rho(x))},~\partial_2 L=\frac{\partial
L}{\partial(\nabla_q^h y(x))},~\partial_3 L= \frac{\partial
L}{\partial(y^\rho(\rho^{\alpha_0}(x)))}~\texttt{and}~\partial_4 L=\frac{\partial
L}{\partial(\nabla_q^hy(\rho^{\alpha_0}(x)))}.$$ and where Lemma
\ref{th} and that $\nabla_q^h \rho^{\alpha_0}(t)=q^{\alpha_0}$ are
used. If we use the change of variable $u=\rho^{\alpha_0}(x)$,
which is a linear function, and make use of Theorem 1.98 in
\cite{Martinbook} and Lemma \ref{th} we then obtain
$$\delta J(y(x),\eta(x))=\int_a^b [\partial_1 L(x)
\eta^\rho(x)+\partial_2 L(x)\nabla_q^h\eta(x)]\nabla_q^hx +$$
\begin{equation} \label{cal22}
 \int_a^{\rho^{\alpha_0}(b)}[q^{-\alpha_0}\partial_3 L((\rho^{\alpha_0})^{-1}(x))\eta^\rho(x)+
q^{-\alpha_0} \partial_4 L((\rho^{\alpha_0})^{-1}(x))\nabla_q^h \eta
(x)]\nabla_q^hx
\end{equation}
where we have used the fact that $\eta\equiv 0$ on
$[\rho^{\alpha_0}(a),a].$

Splitting the first integral in (\ref{cal22})  and rearranging
will lead to

$$\delta J(y(x),\eta(x))=\int_a^{\rho^{\alpha_0}(b)} [\partial_1 L(x) \eta^\rho(x)$$
$$
+\partial_2 L(x)\nabla_q^h\eta(x)+q^{-\alpha_0}\partial_3 L((\rho^{\alpha_0})^{-1}(x))\eta^\rho(x)+
q^{-\alpha_0} \partial_4 L((\rho^{\alpha_0})^{-1}(x))\nabla_q^h \eta
(x)]\nabla_q^hx $$

\begin{equation} \label{cal3}
+\int_{\rho^{\alpha_0}(b)}^b[\partial_1 L(x)
\eta^\rho(x)+\partial_2 L(x)\nabla_q^h\eta(x)]\nabla_q^hx
\end{equation}
If we make use of part 3 of Theorem \ref{t}  then we reach
$$\delta J(y(x),\eta(x))=$$

$$\int_a^{\rho^{\alpha_0}(b)} \left \{  \partial_2 L(x)\nabla_q^h\eta(x)+
q^{-\alpha_0} \partial_4 L((\rho^{\alpha_0})^{-1}(x))\nabla_q^h \eta
(x)+\nabla_q^h[\int_a^x \partial_1 L(z)\nabla_q^h z.\eta(x)]\right.$$
$$-\int_a^x \partial_1 L(z)\nabla_q^h z
.\nabla_q^h\eta(x)+q^{-\alpha_0}\nabla_q^h[\int_a^x
\partial_3 L((\rho^{\alpha_0})^{-1}(z))\nabla_q^h z.\eta(x)]- $$

\begin{equation} \label{off}
\left. q^{-\alpha_0}\int_a^x
\partial_3 L((\rho^{\alpha_0})^{-1}(z))\nabla_q^h z. \nabla_q^h\eta(x)
\right \}\nabla_q^hx +\int_{\rho^{\alpha_0}(b)}^b
\end{equation}
\begin{equation} \label{offf}
\left \{  \partial_2 L(x)\nabla_q^h\eta(x)
+\nabla_q^h[\int_{\rho^{\alpha_0}(b)}^x\partial_1 L(z)\nabla_q^h
z.\eta(x)]\right.-\left.\int_{\rho^{\alpha_0}(b)}^x
\partial_1 L(z)\nabla_q^h z.\nabla_q^h\eta(x)\right \}\nabla_q^hx
\end{equation}
In the above equations (\ref{off}),(\ref{offf}), once choose
$\eta$ such that $\eta(a)=0$ and $\eta \equiv 0$ on
$[q^{\alpha_0}b,b]$ and in another case choose $\eta$ such that
$\eta(b)=0$ and $\eta \equiv 0$ on $[a,q^{\alpha_0}b]$, and then
make use of Lemma \ref{ftcv} to arrive at the following theorem:

\begin{theorem} \label{main1}
Let   $J:S\rightarrow \mathbb{R}$ be the $\mathbb{T}_q^h$-integral
functional
\begin{equation}\label{1main1}
J(y)=\int_a^b L(x,y^\rho(x),\nabla_q
^hy(x),y^\rho(\rho^{\alpha_0}(x)),\nabla_q^hy((\rho^{\alpha_0}(x))\nabla_q^hx
\end{equation}
where
$$a,b\in \mathbb{T}_q^h,a < \rho^{\alpha_0}(b)< b $$

$$L:[a,b]\times (\mathbb{R}^n)^4\rightarrow \mathbb{R}~\mbox{and}
~y^\rho(x)=y(\rho(x))$$

and

$$S=\{y:[\rho^{\alpha_0}(a),b]\rightarrow
\mathbb{R}^n:y(x)=\varphi(x)~(\forall x \in [\rho^{\alpha_0}(a),
a])~\mbox{and}~y(b)=c_0\}.$$ Then the necessary condition for
$J(y)$ to possess an extremum for a given function $y(x)$ is that
$y(x)$ satisfies the following Euler-Lagrange equations
\begin{equation}\label{maine1}
\nabla_q^h\partial_2 L(x)+ q^{-\alpha_0}\nabla_q^h
\partial_4 L((\rho^{\alpha_0})^{-1}(x))=\partial_1 L(x)+q^{-\alpha_0}\partial_3 L((\rho^{\alpha_0})^{-1}(x)),
~(x \in ~[a,\rho^{\alpha_0}(b)]_\kappa)
\end{equation}
and
\begin{equation}\label{maine2}
\nabla_q^h \partial_2 L(x)=\partial_1 L(x)~(x \in ~[\rho^{\alpha_0}(b),b]_\kappa)
\end{equation}
Furthermore, the equation
\begin{equation}\label{maine3}
q^{-\alpha_0}\partial_4 L((\rho^{\alpha_0})^{-1}(x))\eta(x)|_a^{\rho^{\alpha_0}(b)}=0
\end{equation}
holds along $y(x)$ for all admissible variations $\eta(x)$
satisfying $\eta(x)=0, ~~x\in [\rho^{\alpha_0}(a), a]\cup \{b\}$.
\end{theorem}
The necessary condition represented by (\ref{maine3}) is obtained
by applying integration by parts in (\ref{cal3}) and then
substituting equations (\ref{maine1}) and (\ref{maine2}) in the
resulting integrals. The above theorem can be generalized as
follows:
\begin{theorem} \label{maing}
Let   $J:S^m\rightarrow \mathbb{R}$ be the
$\mathbb{T}_q^h$-integral functional

$$J(y_1,y_2,...,y_m)=\int_a^bL(x,y_1^\rho(x),$$
$$ y_2^\rho(x),...,y_m^\rho(x),\nabla_q^hy_1(x),\nabla_q^hy_2(x),
...,\nabla_q^hy_m(x),y_1^\rho(\rho^{\alpha_0}(x)),y_2^\rho(\rho^{\alpha_0}(x))$$
\begin{equation}\label{1main1}
,...,y_m^\rho(\rho^{\alpha_0}(x)),
\nabla_q^hy_1((\rho^{\alpha_0}(x),
\nabla_q^hy_1((\rho^{\alpha_0}(x),...,\nabla_q^hy_m((\rho^{\alpha_0}(x))\nabla_q^hx
\end{equation}
where
$$a,b\in \mathbb{T}_q^h,a < \rho^{\alpha_0}(b)< b $$

$$L:[a,b]\times (\mathbb{R}^n)^{4m}\rightarrow \mathbb{R}~\texttt{and}
~y^\rho(x)=y(\rho(x))$$

and
$$S^m=\{y=(y_1,y_2,...,y_m):$$

$$y_i:[\rho^{\alpha_0}(a),b]\rightarrow
\mathbb{R}^n,y_i(x)=\varphi_i(x)~(\forall x \in
[\rho^{\alpha_0}(a), a])~\texttt{and}~y_i(b)=c_i,i=1,2,...,m\}.$$
Then a necessary condition for $J(y)$ to possess an extremum for
a given function $y(x)=(y_1(x),y_2(x),...,y_m(x))$ is that  $y(x)$
satisfies the following Euler-Lagrange equations
$$\nabla_q^h \partial_2 L^i(x)+ q^{-\alpha_0}\nabla_q^h
\partial_4 L^i((\rho^{\alpha_0})^{-1}(x))=$$
\begin{equation}\label{gmaine1}
\partial_1 L^i(x)+q^{-\alpha_0}\partial_3 L^i((\rho^{\alpha_0})^{-1}(x)),
~(x \in ~[a,\rho^{\alpha_0}(b)]_\kappa)
\end{equation}
and
\begin{equation}\label{gmaine2}
\nabla_q^h \partial_2 L^i(x)=\partial_1 L^i(x)~(x \in
[\rho^{\alpha_0}(b),b]_\kappa)
\end{equation}
Furthermore, the equations
\begin{equation}\label{gmaine3}
q^{-\alpha_0}\partial_4 L^i((\rho^{\alpha_0})^{-1}(x))\eta_i(x)|_a^{\rho^{\alpha_0}(b)}=0
\end{equation}
hold along $y(x)$ for all admissible variations $\eta_i(x)$
satisfying
$$\eta_i(x)=0, ~~x\in [\rho^{\alpha_0}(a), a]\cup \{b\}, i=1,2,...,m$$

where

$$\partial_1 L^i=\frac{\partial L}{\partial(y_i^\rho(x))},~\partial_2 L^i=\frac{\partial
L}{\partial(\nabla_q^h y_i(x))},~\partial_3 L^i= \frac{\partial
L}{\partial(y_i^\rho(\rho^{\alpha_0}(x)))}~\texttt{and}~\partial_4 L^i=\frac{\partial
L}{\partial(\nabla_q ^hy_i(\rho^{\alpha_0}(x)))}.$$
\end{theorem}
\section{The optimal-control problem} \label{s:4}
Our aim in this section is to find the optimal control variable
$u(x)$ defined on the $H-$time scale, which minimizes the
performance index
\begin{equation}\label{optim}
J(y,u)=\int_a^b
L(x,y^\rho(x),u^\rho(x),y^\rho(\rho^{\alpha_0}(x)),\nabla_q^hy((\rho^{\alpha_0}(x))\nabla_q^hx
\end{equation}
subject to the constraint
\begin{equation}\label{constr}
\nabla_q^hy(x)=G(x,y^\rho(x),u^\rho(x))
\end{equation}
such that
\begin{equation}\label{delconst}
y(b)=c,~~y(x)=\phi(x)~~x \in [\rho^{\alpha_0}(a), a])
\end{equation}
$$a,b\in \mathbb{T}_q^h,a < \rho^{\alpha_0}(b)< b $$

$$L:[a,b]\times (\mathbb{R}^n)^4\rightarrow \mathbb{R}~\mbox{and}
~y^\rho(x)=y(\rho(x))$$ where $c$ is a constant and $L$ and $G$
are functions with continuous first and second partial derivatives
with respect to all of their arguments. To find the optimal
control, we define a modified performance index as
$$I(y,u)=\int_a^b
[L(x,y^\rho(x),u^\rho(x),y^\rho(\rho^{\alpha_0}(x)),\nabla_q^hy((\rho^{\alpha_0}(x))$$
\begin{equation} \label{perf}
+\lambda^\rho(x)
(\nabla_q^hy(x)-G(x,y^\rho(x),u^\rho(x)))]\nabla_q^hx
\end{equation}
where $\lambda$ is a Lagrange multiplier or an adjoint variable.

 Using the the equations (\ref{maine1}), (\ref{maine2}) and
(\ref{maine3}) of Theorem \ref{maing} with
$m=3,~~(y_1=y,~y_2=u,~y_3=\lambda)$, the necessary conditions for
our optimal control are (we remark that as there is no any time
scale derivative of $u(x)$, no boundary constraints for it are
needed)

$$ \nabla_q^h \lambda^\rho(x)+q^{-\alpha_0}\nabla_q^h \frac{\partial
L}{\partial\nabla_q^h(y(((\rho^{\alpha_0}(x))}(\rho^{\alpha_0})^{-1}(x))+\lambda^\rho
(x) \frac{\partial G}{\partial y^\rho(x)}-\frac{\partial
L}{\partial y^\rho (x)}$$
\begin{equation} \label{op1}
-q^{-\alpha_0}\frac{\partial
L}{\partial(y^\rho((\rho^{\alpha_0}(x))}((\rho^{\alpha_0})^{-1}(x))=0,~~(x
\in [a,\rho^{\alpha_0}(b)]_\kappa),
\end{equation}

\begin{equation} \label{op2}
\nabla_q^h \lambda^\rho(x)+\lambda^\rho(x) \frac{\partial
G}{\partial y^\rho(x)}-\frac{\partial L}{\partial y^\rho(x)}=0,
~~(x \in [\rho^{\alpha_0}(b),b]_\kappa)
\end{equation}
 \begin{equation} \label{op3}
 \lambda^\rho(x) \frac{\partial G}{\partial u^\rho(x)}-\frac{\partial
L}{\partial u^\rho(x)}=0,~~(x \in [a,b])
\end{equation}
and
\begin{equation} \label{op4}
\frac{\partial
L}{\partial\nabla_q^h(y(((\rho^{\alpha_0}(x))}(\rho^{\alpha_0})^{-1}(x)\eta(x)|_a^{\rho
^{\alpha_0}(b)}=0
\end{equation}

and also

$$\nabla_q^hy(x)=G(x,y^\rho(x),u^\rho(x))$$

Note that the condition (\ref{op4}) disappears when the Lagrangian
$L$ is free of the delayed time scale derivative of $y$.

\section{The discrete and quantum cases} \label{s:5}

We recall that the results in the previous sections are valid for
time scales whose backward jump operator $\rho$ has the form
$\rho(x)=qx-h$, in particular for the time scale $\mathbb{T}_q^h$.

(i) \emph{The discrete case}: If $q=1$ and $h>0$ (of special
interest the case when $h=1$), then our work becomes on the
discrete time scale $h\mathbb{Z}=\{hn:n \in \mathbb{Z}\}$. In this
case the functional under optimization will have the form
$$J^h(y)= h \sum_{i=a+1}^b L(ih,y((i-1)h),\nabla^h
y(ih),y(ih-(d+1)h),\nabla^h y(ih-dh))$$
$$a,b \in \mathbb{Z},d \in \mathbb{N}~\mbox{and}~a<b-d<b,$$  and that
$y(bh)=c,~y(ih)=\varphi(ih)~\mbox{for}~a-d\leq i \leq a$

where
$$\nabla^h y(x)=y(x)-y(x-h),~~x \in h\mathbb{Z}.$$

The necessary condition for $J^h(y)$ to possess an extremum for a
given function $y:\{ih:i=a-d,a-d+1,...,a,a+1,...,b\}\rightarrow
\mathbb{R}^n$ is that  $y(x)$ satisfies the following
h-Euler-Lagrange equations
\begin{equation}\label{dis1}
\nabla^h \partial_2 L(ih)+ \nabla^h \partial_4 L((i+d)h)=\partial_1 L(ih)+\partial_3 L((i+d)h),
~(i=a+1,a+2,...,b-d)
\end{equation}
and
\begin{equation}\label{dis2}
\nabla^h \partial_2 L(ih)=\partial_1 L(ih)~(i=b-d+1,b-d+2,..,b)
\end{equation}
Furthermore, the equation
\begin{equation}\label{dis3}
\partial_4 L(bh)\eta((b-d)h)-\partial_4 L((a+d)h)\eta (ah)=0
\end{equation}
holds along $y(x)$ for all admissible variations $\eta(x)$
satisfying $\eta(ih)=0, ~~i\in \{a-d,a-d+1,..., a\}\cup \{b\}$.

In this case the h-optimal-control problem would read as:

Find the optimal control variable $u(x)$ defined on the time scale
$h\mathbb{Z}$, which minimizes the h-performance index
\begin{equation}\label{hoptim}
J^h(y,u)=h \sum_{i=a+1}^b
L(ih,y((i-1)h),u((i-1)h),y(ih-(d+1)h),\nabla^h y(ih-dh))$$
$$a,b \in \mathbb{Z},d \in \mathbb{N}~\mbox{and}~a<b-d<b,
\end{equation}
subject to the constraint
\begin{equation}\label{hconstr}
\nabla^hy(ih)=G(ih,y((i-1)h),u((i-1)h)),~~i=a+1,a+2,...,b
\end{equation}
such that
\begin{equation}\label{delconst}
y(bh)=c,~~y(ih)=\phi(ih),~~(i=a-d,a-d+1,...,a)
\end{equation}
$$a,b \in \mathbb{N},a < b-d< b $$

 The necessary conditions for this h-optimal control are

$$ \nabla^h \lambda((i-1)h)+\nabla^h \frac{\partial
L}{\partial\nabla^hy((i-d-1)h)}((i+d)h)+\lambda ((i-1)h)
\frac{\partial G}{\partial y((i-1)h)}-$$
\begin{equation} \label{hop1}
\frac{\partial L}{\partial y ((i-1)h)}-\frac{\partial
L}{\partial(y((i-d-1)h)}((i+d)h)=0,~~(i=a+1,a+2,...,b-d),
\end{equation}

\begin{equation} \label{hop2}
\nabla_q^h \lambda((i-1)h)+\lambda((i-1)h) \frac{\partial
G}{\partial y((i-1)h)}-\frac{\partial L}{\partial y((i-1)h)}=0,
~~(i=b-d+1,b-d+2,...,b)
\end{equation}
 \begin{equation} \label{op3}
 \lambda((i-1)h) \frac{\partial G}{\partial u((i-1)h)}-\frac{\partial
L}{\partial u((i-1)h)}=0,~~(i=a,a+1,...,b)
\end{equation}
and
\begin{equation} \label{hop4}
\frac{\partial
L}{\partial\nabla^hy((i+d)h)}(bh)\eta((b-d)h)-\frac{\partial
L}{\partial\nabla^hy((i+d)h)}((a+d)h)\eta(ah)=0
\end{equation}
and also

$$\nabla^hy(ih)=G(ih,y((i-1)h),u((i-1)h)),~~i=a+1,a+2,...,b$$

Note that the condition (\ref{hop4}) disappears when the
Lagrangian $L$ is independent of the delayed $\nabla^h$ derivative
of $y$.

\begin{example} \label{disc}
In order to illustrate our results we analyze an example of
physical interest. Namely, let us consider the following discrete
action,
$$J^h(t)=\frac{h}{2}\sum_{i=a+1}^b [\nabla^h y(ih)]^2-V(y(ih-(d+1)h)),
~a,b \in \mathbb{N},~a<b-d<b$$ subject to the condition
$$y(bh)=c,~~y(ih)=\varphi(ih),~~\mbox{for}~~i=a-d,a-d+1,...,a$$
The corresponding h-Euler-Lagrange equations are as follows:
\begin{equation} \label{disc1}
y(ih)-2y((i-1)h)+y((i-2)h)+\frac{\partial V}{\partial
y(ih-(d+1)h)}((i+d)h),~~(i=a+1,...,b-d)
\end{equation}
and
\begin{equation} \label{disc2}
y(ih)-2y((i-1)h)+y((i-2)h)=0, ~~(i=b-d+1,b-d+2,...,b)
\end{equation}
We observe that when the delay is removed, that is $d=0$, the
classical discrete Euler-Lagrange equations are reobtained.
\end{example}

(ii) \emph{The quantum case}: If $0<q<1$ and $h=0$, then our work
becomes on the  time scale  $\mathbb{T}_q=\{q^n:n \in
\mathbb{Z}\}\cup \{0\}$. In this case the functional under
optimization will have the form

\begin{equation}\label{qm1}
J_q(y)=\int_a^b L(x,y(qx),\nabla_q
y(x),y(xq^{\alpha_0+1}),\nabla_qy(xq^{\alpha_0}))\nabla_qx
\end{equation}

where
$$a=q^{\alpha+1},b=q^{\beta}, \alpha,\beta, \alpha_0 \in \mathbb{Z},
\alpha > \beta~ \mbox{and} ~~ \beta+\alpha_0 < \alpha,$$

$$L:[a,b]_q\times (\mathbb{R}^n)^4\rightarrow
\mathbb{R}~\mbox{and}~[a,b]_q=\{q^i:i=\alpha+1,\alpha+2,...,\beta\}
~$$

Using the $\nabla-$integral theory on time scales,  the functional
$J_q$ in (\ref{qm1}) turns to be

$$J_q(y)=(1-q)\sum_{i=\alpha}^\beta q^i L(q^i,y(q^{i+1}),\nabla_q
y(q^i),y(q^{\alpha_0 +i+1}),\nabla_q y(q^{\alpha_0 +i}) )$$

The necessary condition for $J_q(y)$ to possess an extremum for a
given function
$y:\{q^i:i=\alpha+1-\alpha_0,\alpha+2-\alpha_0...,\alpha,\alpha+1,...,\beta\}\rightarrow
\mathbb{R}^n$ is that  $y(x)$ satisfies the following
q-Euler-Lagrange equations

$$\nabla_q\frac{\partial L}{\partial(\nabla_q y(x))}(x)+
q^{-\alpha_0}\nabla_q \frac{\partial
L}{\partial(\nabla_qy(q^{\alpha_0}x))}(q^{-\alpha_0}x)=$$

\begin{equation}\label{quantum1}
\frac{\partial L}{\partial y(qx)}(x)+q^{-\alpha_0}\frac{\partial
L}{\partial y(q^{\alpha_0+1}x)}(q^{-\alpha_0}x), ~(x \in
[a,q^{\alpha_0}b]_\kappa)
\end{equation}
and
\begin{equation}\label{quantum2}
\nabla_q\frac{\partial L}{\partial(\nabla_q
y(x))}(x)=\frac{\partial L}{\partial(y(qx))}(x)~(x \in
[q^{\alpha_0}b,b]_\kappa)
\end{equation}
Furthermore, the equation
\begin{equation}\label{quantum3}
q^{-\alpha_0}\frac{\partial
L}{\partial(\nabla_qy(\rho^{\alpha_0}(x)))}(q^{-\alpha_0}x)\eta(x)|_a^{q^{\alpha_0}b}=0
\end{equation}
holds along $y(x)$ for all admissible variations $\eta(x)$
satisfying $\eta(x)=0, ~~x\in [q^{\alpha_0}a, a]_q\cup \{b\}$.

In this case the q-optimal-control problem would read as:

Find the optimal control variable $u(x)$ defined on the
$\mathbb{T}_q-$time scale, which minimizes the performance index
\begin{equation}\label{qoptim}
J_q(y,u)=\int_a^b
L(x,y(qx),u(qx),y(q^{\alpha_0+1}x),\nabla_qy(q^{\alpha_0+1}x)\nabla_qx
\end{equation}
subject to the constraint
\begin{equation}\label{qconstr}
\nabla_q^hy(x)=G(x,y(qx),u(qx))
\end{equation}
such that
\begin{equation}\label{qdelconst}
y(b)=c,~~y(x)=\phi(x)~~x \in [q^{\alpha_0}a, a])
\end{equation}
$$a=q^{\alpha+1},b=q^\beta, \alpha_0+\beta < \alpha +1$$

$$L:[a,b]_q\times (\mathbb{R}^n)^4\rightarrow \mathbb{R}$$
where $c$ is a constant and $L$ and $G$ are functions with
continuous first and second partial derivatives with respect to
all of their arguments.

 The necessary conditions for this q-optimal control are:

$$ \nabla_q \lambda(qx)+q^{-\alpha_0}\nabla_q\frac{\partial
L}{\partial\nabla_qy(q^{\alpha_0+1}x)}(q^{-\alpha_0}x)+\lambda
(qx) \frac{\partial G}{\partial y(qx)}-\frac{\partial L}{\partial
y(qx)}$$
\begin{equation} \label{qop1}
-q^{-\alpha_0}\frac{\partial L}{\partial
y(q^{\alpha_0+1}x)}(q^{-\alpha_0}x)=0,~~(x \in
[a,q^{\alpha_0}b]_\kappa),
\end{equation}

\begin{equation} \label{qop2}
\nabla_q \lambda(qx)+\lambda(qx) \frac{\partial G}{\partial
y(qx)}-\frac{\partial L}{\partial y^\rho(x)}=0, ~~(x \in
[q^{\alpha_0}b,b]_\kappa)
\end{equation}
 \begin{equation} \label{qop3}
 \lambda(qx) \frac{\partial G}{\partial u(qx)}-\frac{\partial L}{\partial
u(qx)}=0,~~(x \in [a,b]_q)
\end{equation}
and
\begin{equation} \label{qop4}
\frac{\partial
L}{\partial\nabla_qy(q^{\alpha_0}x)}(q^{-\alpha_0}x)\eta(x)|_a^{q
^{\alpha_0}b}=0
\end{equation}

and also

$$\nabla_q^hy(x)=G(x,y(qx),u(qx))$$

Note that the condition (\ref{qop4}) disappears when the
Lagrangian $L$ is independent of the delayed $\nabla_q$ derivative
of $y$.

\begin{example} \label{exquantum}
Suppose that the problem is that of finding a control function
$u(x)$ defined on the time scale $T_q$ such that the corresponding
solution of the controlled system
\begin{equation} \label{qcontrol}
\nabla_q y(x)= -ry(qx)+u(qx),~~~r>0,
\end{equation}
satisfying the conditions:
$$y(b)=c,~~y(x)=\varphi(x),~~\mbox{for}~x \in
[q^{\alpha_0}a,a]_q,~~a=q^{\alpha+1},b=q^\beta, \alpha_0+\beta <
\alpha +1$$ is an extremum for the q-integral functional
(q-quadratic delay cost functional):
\begin{equation} \label{qcontrol2}
J_q(y(x),u(x))=\frac{1}{2}(1-q)\sum_{i=\alpha}^\beta q^i
[y^2(q^{i+\alpha_0+1})+u^2(q^{i+1})]
\end{equation}
According to (\ref{qop3}) and (\ref{qop4}), the solution of the
problem satisfies:
\begin{equation} \label{sat1}
\nabla_q \lambda (qx)=r \lambda (qx)+q^{-\alpha_0}y(qx),~~~(x \in
[a,q^{\alpha_0}b]_\kappa),
\end{equation}
\begin{equation} \label{sat2}
\nabla_q \lambda (qx)=r \lambda (qx),~~(x \in
[q^{\alpha_0}b,b]_\kappa),
\end{equation}
\begin{equation} \label{sat2}
\lambda (qx)=u(qx),~~~( x \in [a,b]_q)
\end{equation}
and of course
$$\nabla_q y(x)= -ry(qx)+u(qx)$$
\end{example}
When the delay is absent (i.e $\alpha_0=0$), it can be shown that
the above system is reduced to a second order q-difference
equation. Namely, reduced to
 $$\nabla_q^2 y(x)+rq(\nabla_q y)(qx)=q (r^2+1)y(qx)+qr \nabla_q y(x)$$
 If we solve recursively for this equation in terms of an integer power
series by using the initial data, then the resulting solution will
tend to the solutions of the second order linear differential
equation:
$$y^{\prime\prime}-(r^2+1)y=0.$$
Clearly the solutions for this equation are : $exp(\sqrt{r^2+1}x)$
and $exp(-\sqrt{r^2+1}x)$. For details see \cite{Gas}.

\section{Conclusion} \label{s:6}

In this manuscript we have developed an optimal variational
problem in the presence of delay on time scales whose backward
jumping operators are of the form $\rho(t)=qt-h,~q>0,~h\geq 0$,
called $H-$time scales. Such kinds of time scales unify the
discrete, the quantum and the continuous cases, and hence the obtained
results generalized many previously obtained results either in the
presence of delay or without. To formulate the necessary
conditions for this optimal control problem, we first obtained the
Euler-Lagrange equations for one unknown function then generalized
to the n-dimensional case. The state variables of the Lagrangian
in this case  are defined on the $H-$time scale and
contain some delays. When $q=1$ and $h=0$ with the existence of
delay some of the results in \cite{Bliss} are recovered. When
$0<q<1$ and $h=0$ and the  delay is absent most of the results in
\cite{Gas} can be reobtained. When $q=1$ and the delay is absent
some of the results in \cite{Cad} are reobtained. When the delay
is absent and the time scale is free somehow, some of the results
in \cite{Nat} can be recovered as well.

Finally, we would like to mention that we followed the line of
nabla time scale derivatives in this article, analogous results
can be originated if the delta time scale derivative approach is
followed.

\section*{Acknowledgments}
This work is partially supported by the Scientific and Technical
Research Council of Turkey.

  \end{document}